\documentclass[12pt, amscd]{amsart}
\usepackage{graphicx}
\usepackage{amsmath,amssymb, amsthm}
\newtheorem{thm}{Theorem}[section]

\newtheorem{lem}[thm]{Lemma}

\theoremstyle{definition}
\newtheorem{defn}[thm]{Definition}
\newtheorem{exm}[thm]{Example}

\DeclareMathOperator{\girth}{girth}

\DeclareMathOperator{\bht}{bight}
\DeclareMathOperator{\hgt}{ht}

\DeclareMathOperator{\ass}{Ass}
\DeclareMathOperator{\pc}{\mathcal{PC}}

\def\Z {\mathbb Z}

\begin{document}

\title[Cohen-Macaulay edge-weighted graphs of girth $5$ or greater] {Cohen-Macaulay edge-weighted graphs of girth $5$ or greater}
\author[T.T. Hien]{Truong Thi Hien}
\address{Faculty of Natural Sciences, Hong Duc University,
No. 565 Quang Trung Street, Dong Ve Ward, Thanh Hoa, Vietnam}
\email{hientruong86@gmail.com}

\date{}
\commby{}
%-----------------------------------------------------------
\begin{abstract}
Let $G_\omega$ be an edge-weighted graph whose underlying graph is $G$. In this paper, we  enlarge the class of Cohen-Macaulay edge-weighted graphs $G_\omega$ by classifying completely them when the graph $G$ has girth $5$ or greater.
\end{abstract}
% -----------------------------------------------------------
\maketitle
% -----------------------------------------------------------
\section*{Introduction}
Let $R = K[x_1,\ldots,x_d]$ be a standard graded polynomial ring over a given field $K$.  Let $G$ be a simple graph with the vertex set $V(G) = \{x_1,\ldots,x_d\}$ and the edge set $E(G)$.  By abuse of notation, we also use $x_ix_j$ to denote an edge $\{x_i, x_j\}$ of $G$. A \textit{edge-weighted graph} $G_{\omega}$ (whose underlying graph is $G$) is the couple $(G,\omega)$, where $\omega$ is a function $\omega \colon E(G) \to \Z_{>0}$, which is called a {\it weight edge} on $G$. An edge-weighted graph $G_{\omega}$ where each edge has the same weight is a trivial edge-weighted graph. The {\it weighted edge ideal} of $G_\omega$ was introduced by Paulsen and Sather-Wagstaff \cite{PS}, given by 
$$I(G_\omega) = ((x_ix_j)^{\omega(x_ix_j)}\mid x_i x_j\in E(G)).$$

We say that the edge-weighted graph $G_\omega$ was called {\it Cohen-Macaulay}  if $R/I(G_\omega)$ is Cohen-Macaulay. In \cite{PS}, the authors constructed the irreducible decomposition of $I(G_\omega)$ and classified  Cohen-Macalay edge-weighted graphs $G_\omega$ where the underlying graph $G$ is a tree or a cycle.  After that Fakhari, Shibata, Terai and Yassemi classified Cohen-Macalay edge-weighted graphs $G_\omega$ when $G$ is a very well-covered graph (see \cite{SSTY}).  It is worth mentioning that the problem of classifying sequentially Cohen-Macaulay edge-weighted graphs is studied in \cite{MDV}, and classifying Cohen-Macaulay vertex-weighted oriented is studied in \cite{DT, HLMRV, PRT, PRV}. In this paper, we study Cohen-Macaulay properties for the edge-weighted graphs $G_{\omega}$. More specifically, we classify Cohen-Macaulay edge-weighted graphs $G_{\omega}$ when $G$ has girth at least $5$. Recall that the {\it girth} of a graph $G$, denoted by $\girth (G)$, is the length of the shortest cycle contained in it. If a graph contains no cycle, its girth is defined to be infinite. 

The main result of the paper is the following theorem.

\medskip

{\noindent {\bf Theorem \ref{main-theorem}}.  \it Let $G$ be a graph of girth at least $5$ and $\omega$ is a weight edge on $G$. Then, the following conditions are equivalent:
	\begin{enumerate}
		\item $G_\omega$ is Cohen-Macaulay.
		\item $G$ is Cohen-Macaulay and $G_\omega$ is unmixed.
		\item $G$ is in the class $\pc$ and the weight edge $\omega$ on $G$ satisfies:
		\begin{enumerate}
			\item The weight of any pendant edge in $G$ is greater than or equal to the weight of every edge adjacent to it.
			\item Every basic $5$-cycle $C$ of $G$ has a balanced vertex adjacent to two vertices on $C$ of degree $2$.
			\item If a vertex $x$ is on a basic $5$-cycle $C$ with $\deg_G(x)\geqslant 3$ and $N_C(x) = \{y,v\}$, then $\min\{\omega(xy),\omega(xv)\} \geqslant \max\{\omega(xw) \mid w \in N_G(x)\setminus\{y,v\}\}$.
		\end{enumerate}
	\end{enumerate}
}

\medskip

To understand the above theorem clearly, we first recall some definitions and terminologies. An edge-weighted graph  $G_\omega$ is called \textit{unmixed} if the quotient ring $R/I(G_\omega) $ is unmixed. 
 An edge of $G$ is called the {\it pendant edge} if one of its vertices is a leaf. A \textit{basic $5$-cycle} is a cycle of length $5$ and there are no two adjacent vertices of degree three or more in $G$.  

For a given graph $G$, let  $C(G)$ and $P(G)$ denote  the set of all vertices that belong to basic 5-cycles and pendant edges, respectively. $G$ is said to be \textit{in the class $\mathcal{PC}$} if 
\begin{enumerate}
	\item $V(G)$ can be partitioned into $V(G)=P(G) \cup C(G)$; and
	\item the pendant edges form a perfect matching of $G[P(G)]$.  
\end{enumerate}

Let $C$ be an induced $5$-cycle of $G$ with $E(C) = \{xy,yz,zu,uv,vx\}$.  We say that the vertex $x$ is a \textit{balanced vertex} on $C$ (with respect to $\omega$)  if
\begin{enumerate}
	\item $\omega(xy) = \omega(xv)$; and
	\item $\omega(xy) \leqslant \omega(yz)\geqslant \omega(zu) \leqslant \omega(uv) \geqslant \omega(xv)$.
\end{enumerate}

\begin{figure}[h]
	\includegraphics[scale=0.45]{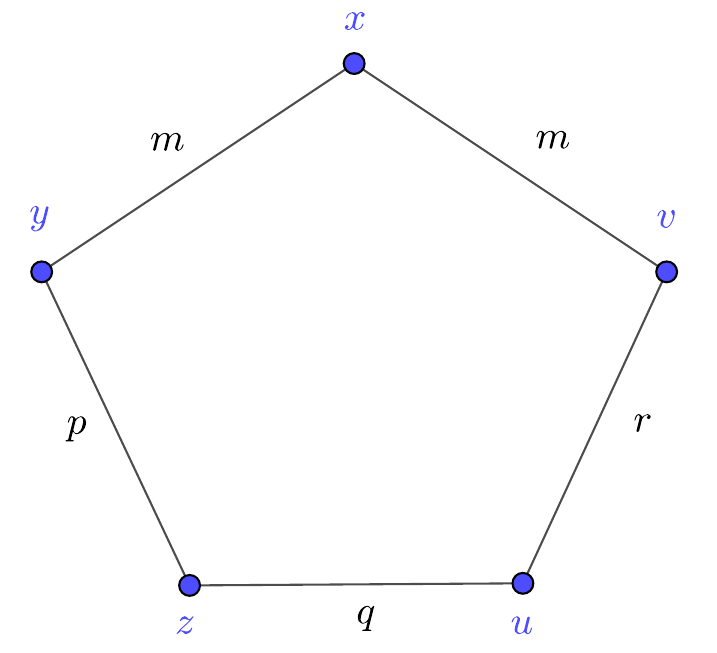}\\
	\medskip
	\caption{The balanced vertex $x$ on $C$.}
	\label{PSW}
\end{figure}

This definition is motivated by \cite[Theorem 4.4]{PS}, which says that $C_\omega$ is Cohen-Macaulay if and only if $C$ has a balanced vertex. In Figure \ref{PSW}, where the weight edge is indicated on edges, $x$ is a balanced vertex on $C$ if the following inequalities hold:  $m \leqslant p \geqslant q \leqslant r\geqslant m$.

\medskip

Let us explain the ideal to prove the theorem \ref{main-theorem}. We will prove this theorem by the following sequence: $(1) \Rightarrow (2) \Rightarrow (3) \Rightarrow (1)$. By \cite[Theorem 2.6]{HTT}, if $G_\omega$ is Cohen-Macaulay, then $I(G) = \sqrt{I(G_\omega)}$ is also Cohen-Macaulay, thus we get $(1) \Rightarrow (2)$.  To prove $(2) \Rightarrow (3)$, we has the result that $G$ is in the class $\pc$ if $G$ has girth at least $5$. In addition, we introduce the notion of weighted vertex cover with minimal support to characterize the associated primes of $I(G_\omega)$. Together with the structure of $G$, we can prove the combinatorial properties $(a)$-$(c)$. It remains to show that $(3) \Rightarrow (1)$. If $G_\omega$ satisfies the condition $(3)$, we will prove $G_\omega$ is Cohen-Macaulay by induction on the number of basic 5-cycles of $G$. Indeed, assume $x$ is a balanced vertex on some basic $5$-cycle $C$ as indicated in the property $(b)$ and $m=\omega(xy)$ with $xy\in E(C)$. We show that $(I(G_\omega),x^m)$ and $I(G_\omega)\colon x^m$ are the weighted edge ideals of some edge-weighted graphs. Furthermore, these edge-weighted graphs also satisfy the condition $(3)$ and have less the number of basic 5-cycles than $G$, then they are Cohen-Macaulay by induction. Therefore, the conclusion is followed.

 \medskip
The paper consists of two sections.  In Section $1$, we set up some basic notations, terminologies from the graph theory, the irreducible decomposition of the weighted edge ideal of an edge-weighted graph, and Cohen-Macaulay monomial ideals and their colon ideals. In Section $2$, we classify  Cohen-Macaulay edge-weighted graphs of girth at least $5$ by giving some characteristics of the weight $\omega$ on pendant edges and basic $5$-cycles of $G$.

\section{Preliminaries}

We begin this section with some observations from the graph theory. Let $G = (V(G), E(G))$ be a simple graph. Note that two vertices of $G$ are adjacent if they are connected by an edge; two edges of $G$ are adjacent if they share a common vertex.

\begin{defn}
A set of vertices is called a \textit{vertex cover} of $G$ if for every edge, $(u, v) \in E(G)$, either $u$ or $v$ or both are a part of the set. A \textit{minimal vertex cover} is a vertex cover that no its proper subset is still a vertex cover.
\end{defn}

\begin{defn}
 The set of non-adjacent vertices is called an \textit{independent set}. A \textit{maximal independent set} is an independent set that is not contained properly in any other independent set of $G$. An independent set is called \textit{maximum} if it is of the largest cardinality.  %Denoted $\alpha(G)$ to be the cardinality of a maximum independent set in $G$ and is called the independence number of $G$.
\end{defn}

\textbf{Remark.} Obviously, a vertex cover corresponds to the complement of an independent vertex set. 

 \begin{defn}
 	A subset $P$ of edges of $G$ is a \textit{matching} if there are no two edges in $P$ which are adjacent to each other. A matching $P$ of $G$ is \textit{perfect} if  every vertex of $G$ is incident to some edge in $P$, i.e. in the case $|V(G)| = 2|P|$.
 \end{defn}

\medskip

If $X \subseteq V(G)$, $G[X]$ is the induced subgraph of $G$ on $X$. By $G\setminus X$, we mean the induced subgraph $G[V \setminus X]$. The \textit{neighbor} of a vertex $v$ of $G$ means the vertices that are adjacent to $v$ in $G$. The \textit{(open) neighborhood} of a vertex $v$ is the set of its neighbors, i.e., $N_G(v) = \{w \mid w \in V(G) \text{ and } vw\in E(G)\}$. The \textit{closed neighborhood} of $v$ means to all the neighbors of $v$ and itself, i.e., $N_G[v] = N_G(v) \cup \{v\}$; if there is no ambiguity on $G$, we use $N(v)$ and $N[v]$, respectively. We also use the symbol $N_G[X] = X \cup \{v\mid vu \in E(G) \text{ for some } u\in X\}$ to denote the closed neighborhood of $X$ in $G$. The {\it degree} of $v$ in $G$ is the number of its neighbors and is denoted by $\deg_G(v)$. It implies that $\deg_G(v) = |N_G(v)|$. Note that $v$ is called a leaf if $\deg_G(v) = 1$. 

\medskip

We next introduce the class of vertex decomposable graphs (see e.g. \cite{W}). For a vertex $v$ of $G$, denoted $G\setminus v = G\setminus\{v\}$ and $G_v=G\setminus N_G[v]$.

\begin{defn}
 A graph $G$ is called \textit{vertex decomposable}  if it is a totally disconnected graph (i.e. with no edges) or there is a vertex $v$ in $G$ such that 
\begin{enumerate}
	\item $G\setminus v$ and $G_v$ are both vertex decomposable, and
	\item for every independent set $S$ in $G_v$, there is some $u \in N_G(v)$ such that $S \cup \{u\} $ is independent in $G\setminus v$.
\end{enumerate}
\end{defn}
 
 The vertex $v$ which satisfies the condition $(2)$ is called a {\it shedding vertex} of $G$. Recall a graph $G$ is well-covered (see \cite{P}) if every maximal independent set of $G$ has the same size, namely $\alpha(G)$. Thus, if $G$ is well-covered and $v$ is a shedding vertex of $G$, then $G\setminus v$ is also well-covered with $\alpha(G\setminus x) =\alpha(G)$.  
 
 \medskip

Now, we consider some results of the irreducible decomposition of the weighted edge ideal of an edge-weighted graph, and Cohen-Macaulay monomial ideals and their colon ideals, which we shall need in the proof of the main theorem. 

 It is widely known that (see e.g. \cite[Proposition 6.1.16]{Vi})  $$\ass(R/I(G)) = \{(v\mid v \in C) \mid C \text{ is a minimal vertex cover of } G\}.$$ 
Particularly,  $\dim R/I(G) = \alpha(G)$ whenever $V(G)=\{x_1,\ldots,x_d\}$. The graph  $G$ is called a Cohen-Macaulay graph if the ring $R/I(G)$ is Cohen-Macaulay. In consequence, $G$ is well-covered if it is Cohen-Macaulay. 

Let $G_\omega$ be an edge-weighted graph. We know that the usual edge ideal of $G$, denoted by $I(G)$, is a special case of the weighted edge ideal $I(G_\omega)$ when the weight $\omega$ on $G$ is the trivial one, i.e., $\omega(e)=1$ for all $e\in E(G)$. Since $I(G)=\sqrt{I(G_\omega)}$, by \cite[Theorem 2.6]{HTT}, $G$ is Cohen-Macaulay if so is $G_\omega$. Therefore, if we know the structure of the underlying Cohen-Macaulay graph together with the weight edges on it, we can get the picture of the Cohen-Macaulayness of an edge-weighted graph. In this paper, we consider graphs of girth at least $5$ and so the following result plays a crucial role in the paper (see \cite[Theorem 20]{BC} or \cite[Theorem 2.4]{HMT}).

\begin{lem}  \label{HMT} Let $G$ be a connected graph of girth at least $5$. Then, the following statements are equivalent:
\begin{enumerate}
\item $G$ is well covered and vertex decomposable;
\item $G$ is Cohen-Macaulay;
\item $G$ is either a vertex or in the class $\mathcal{PC}$.
\end{enumerate}
\end{lem}

We next describe the associated primes of $R/I(G_\omega)$. 

\begin{defn}
Let $G_\omega$ be an edge-weighted graph, $C$ be a vertex cover of $G$ and a function $\delta \colon C \to \Z_{>0}$. The pair $(C,\delta)$ is called a \textit{weighted vertex cover} of $G_{\omega}$  if for every $e = uv\in E(G)$ we have either $u\in C$ and $\delta(u) \leqslant \omega(e)$ or $v\in C$ and $\delta(v)\leqslant \omega(e)$. 
\end{defn} 

Observe that a pair $(C,\delta)$ where $C\subseteq V(G)$ and $\delta\colon C\to \Z_{>0}$ is a weighted vertex cover of $G_\omega$ if and only if  $P(C,\delta) =(v^{\delta(v)} \mid v\in C) \supseteq I(G_\omega)$. Now, we give a definition of an ordering of weighted vertex covers. 

\begin{defn}
Let $G_\omega$ be an edge-weighted graph. For two weighted vertex covers $(C,\delta)$ and $(C',\delta')$ of $G_\omega$, we say that $(C,  \delta) \leqslant (C',\delta')$ if $C\subseteq C'$ and $\delta(v) \geqslant \delta'(v)$ for every $v\in C$.
\end{defn}

   In the usual sense, $(C,\delta)$ is minimal if it is minimal with respect to this order. Then,

\begin{lem}\label{ID} \cite[Theorem 3.5]{PS} $I(G_\omega)$ can be represented as
$$I(G_\omega)=\bigcap_{(C,\delta) \text{ is minimal }} P(C,\delta)$$
and the intersection is irredundant.
\end{lem}

This lemma  implies that if $(C,\delta)$ is a minimal weighted vertex cover of $G_\omega$, then $(v\mid v\in C) \in \ass(R/I(G_\omega))$.   We say that a weighted vertex cover $(C,\delta)$ of $G_\omega$ is {\it minimal support} if there is no proper subset $C'$ of $C$ such that $(C',\delta)$ is a weighted vertex cover of $G_\omega$.

\begin{lem}\label{cover} If a weighted vertex cover $(C,\delta)$ of $G_\omega$ is minimal support,  then $$(v \mid v \in C) \in \ass(R/I(G_\omega)).$$
\end{lem}
\begin{proof} Since $I(G_\omega)\subseteq P(C,\delta)$,  by Lemma \ref{ID} we have $P(C',\delta') \subseteq P(C,\delta)$ for some minimal weighted vertex cover $(C',\delta')$.  In particular,  $(C',\delta') \leqslant (C,\delta)$.  This implies that $(C', \delta)$ is a weighted vertex cover of $G_\omega$. Since $(C,\delta)$ is minimal support, we have $C= C'$.  Therefore,  $(v\mid v\in C)\in \ass(R/I(G_\omega))$, as required.
\end{proof}

A monomial ideal $I$ is {\it unmixed} if every its associated prime has the same height.  It is well known that if $R/I$ is Cohen-Macaulay, then $I$ is unmixed.  % If $I(G_\omega)$ is unmixed,  we say that $G_\omega$ is unmixed. 
 Because $I(G)=\sqrt{I(G_\omega)}$, hence $G$ is well-covered if $G_\omega$ is unmixed.  In this case,  if $(C,\delta)$ is a weighted minimal vertex cover of $G_\omega$, then $C$ is a minimal vertex cover of $G$.

\medskip

We now recall some techniques to study the Cohen-Macaulayness of monomial ideals as mentioned in \cite{DT}.

\begin{lem}\cite[Lemma 1.4]{DT} \label{CM-Q} Let $I$ be a monomial ideal and $f$ a monomial not in $I$.  We have 
\begin{enumerate}
\item If $I$ is Cohen-Macaulay, then $I\colon f$ is Cohen-Macaulay.
\item If $I\colon f$ and $(I,f)$ are Cohen-Macaulay with $\dim R/I\colon f = \dim R/(I,f)$, then $I$ is Cohen-Macaulay.
\end{enumerate}
\end{lem}

\begin{lem} \cite[Lemma 1.5]{DT} \label{dim} Let $G$ be a well-covered graph.  If $v$ is a shedding vertex of $G$, then 
$$\dim R/I(G) = \dim R/(I(G\setminus v),v)=\dim R/I(G):v.$$
\end{lem}

In the sequel, we need the following lemma obtained from \cite{DT}.

\begin{lem}\label{CM-L02} Let $G$ be a graph in the class $\mathcal{PC}$.  Let $C$ be a basic $5$-cycle and $x$ a vertex in $C$ with degree at least $3$.  Assume that $E(C) = \{xy,yz,zu,uv,vx\}$ and $N(x) = \{y,v,y_1,\ldots,y_k\}$.  Then,  there is an independent set of $G$ with $k$ vertices,  say $\{z_1,\ldots,z_k\}$,  such that
\begin{enumerate}
\item $G[y_1,\ldots,y_k,z_1,\ldots,z_k]$ consists of $k$ disjoint edges $y_1z_1, \ldots, y_kz_k$.
\item $N[z_1,\ldots,z_k] \cap V(C) = \emptyset$.
\end{enumerate}
\end{lem}
\begin{proof} Follows from Part $(1)$ of \cite[Lemma 2.2]{DT}.
\end{proof}

\section{Cohen-Macaulay edge-weighted graphs}

In this section, we classify the Cohen-Macaulay edge-weighted graphs $G_\omega$ of girth at least $5$. Because $G$ is in the class $\pc$, which $V(G) = P(G) \cup C(G)$, then we will study the weight $\omega$ on pendant edges and basic $5$-cycles of $G$ as a natural.

To investigate the weight on pendant edges, we consider the following lemma.

\begin{lem} \label{CM-L01} Let $G_\omega$ be an unmixed edge-weighted graph.   Assume that $(xy)^{\omega(xy)}$ and $(xz)^{\omega(xz)}$ are among  minimal generators of $I(G_\omega) \colon f$ for some monomial $f\notin I(G_\omega)$.  Assume that $x^k \notin I(G_\omega) \colon f$ for every $k$.  If $y$ does not appear in any minimal generator of $I(G_\omega)\colon f$ except for $(xy)^{\omega(xy)}$,   then 
$\omega(xy) \geqslant \omega(xz)$.
\end{lem}
\begin{proof} Follows from \cite[Lemma 2.1]{DT}.
\end{proof}

We now move on to  investigate the weight of basic $5$-cycles.

\begin{lem} \label{CM-L03}  Let $G_\omega$ be an unmixed edge-weighted graph where $G$ is in the class $\pc$. Assume that $C$ is a basic $5$-cycle of $G$ such that that $E(C) = \{xy,yz,zu,uv,vx\}$ and $\deg(x) >2$. Then,
\begin{enumerate}
\item $\omega(xw) \leqslant \min\{\omega(xy), \omega(xv)\}$ for all $w \in N(x) \setminus \{y,v\}$.
\item $\omega(zu) \leqslant \min \{\omega(zy), \omega(uv)\}$.
\end{enumerate}
\end{lem}
\begin{proof} Let $m = \omega(xy),  n = \omega(xv),  p = \omega(yz),  q = \omega(zu),  r = \omega(uv)$.  Assume that
$$N(x) = \{y,v,y_1,\ldots,y_k\}, \text{ where } k \geqslant 1,$$
and $m_i = \omega(xy_i)$ for $i=1,\ldots,k$ so that $m_1\geqslant m_2\geqslant \cdots\geqslant m_k$.  

$(1)$ Assume on the contrary that $m_1 > \min\{m,n\}$. We may assume that $m_1 > m$.  By Lemma \ref{CM-L02},  there is an independent set $\{z_1, \ldots,z_k\}$ of $G$ such that the graph $G[y_1,\ldots,y_k,z_1,\ldots,z_k]$ consists of disjoint edges $y_1z_1,\ldots, y_kz_k$ and $N[z_1,\ldots,z_k] \cap V(C) = \emptyset$.

\begin{figure}[h]
\includegraphics[scale=0.5]{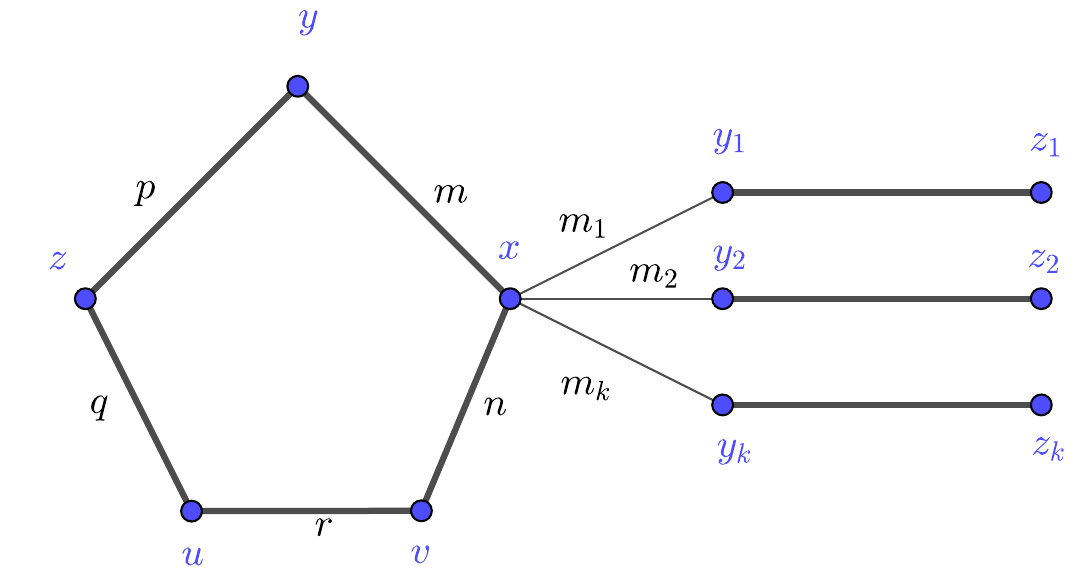}\\
\medskip
\caption{The structure of $G$.}
\label{structureG}
\end{figure}

Let $S_1 = \{z_2,\ldots,z_k\}$.  As $\girth(G) \geqslant 5$, we deduce that $\{y_1, y, u\} \cup S_1$ is an independent set of $G$.  Now extend this set to a maximal independent set of $G$,  say $S$.  Then, $C^* = V(G) \setminus S$ is a minimal cover of $G$.  In particular, $\hgt(I(G_\omega)) = |C^*|$. Let $\delta \colon C^* \to \Z_{>0}$ be such that $(C^*,\delta)$ is a weighted vertex cover of $G_\omega$.  Note that $z,v,x,  y_2, \ldots,y_k \in C^*$ and $y_1, y,u \notin C^*$.   

Let $C' = C^*\cup \{y\}$ and $\delta' \colon C' \to \Z_{>0}$ defined by
$$\delta'(w) = \begin{cases} m_1 & \text{ if } w = x,\\
m & \text{ if } w = y,\\
\min\{n, \delta(v)\} & \text{ if } w = v,\\
\min\{m_i, \delta(y_i)\} & \text{ if } w = y_i,  \text{ for } i = 2,\ldots,k,\\
\delta(w) & \text{ otherwise}.
\end{cases}
$$

We now prove that  $(C', \delta')$ is a weighted vertex cover of $G_\omega$. Indeed, since  $C' = C^*\cup \{y\}$, $(C^*,\delta)$ is a weighted vertex cover of $G_\omega$ and by the definition of the function $\delta'$, it suffices to check the condition of a weighted vertex cover of $G_\omega$ for the set of edges 
$$\{xy, xv, xy_1, xy_2, \ldots, xy_k, yz, vu\}.$$

If $e = xy$, then $y \in C'$ and $\delta'(y) = m = \omega(xy)$.

If $e = xv$, then $v \in C'$ and $\delta'(v) = \min\{n, \delta(v)\} \leqslant n = \omega(xv)$.

If $e = xy_1$, then $x \in C'$ and $\delta'(x) = m_1 = \omega(xy_1)$.

If $e = xy_i$, then $y_i \in C'$ and $\delta'(y_i) = \min\{m_i, \delta(y_i)\} \leqslant m_i = \omega(xy_i)$, for $i = 2, \ldots, k$.

If $e = yz$. By considering the weighted vertex cover $(C^*,\delta)$ we have $\delta(z) \leqslant \omega(yz)$ since $y \notin C^*$ and $z \in C^*$. Thus, for $(C', \delta')$, we have  $z \in C'$ and $\delta'(z) = \delta(z) \leqslant \omega(yz)$.

If $e = vu$, similarly as the previous case, we have $v \in C'$ and $\delta'(v) = \min\{n, \delta(v)\} \leqslant \delta(v) \leqslant \omega(uv)$.

Therefore, $(C', \delta')$ is a weighted vertex cover of $G_\omega$, as desired.

Next, we claim that it is minimal support.  Indeed, assume on the contrary that it is not the case, then there is a vertex, say $w\in C'$ such that $(C'\setminus \{w\},\delta')$ is still a weighted vertex cover of $G_\omega$.  We consider the case $w \in \{y,x,v,z\}$.  Since $u,y_1\notin C'$, it follows that $w$ must be $y$. But in this case,  we have $\delta'(x) = m_1 \leqslant \omega(xy)=m$ if we look at the edge $e= xy$,  a contradiction. Thus, $w\notin  \{y,x,v,z\}$.  Note that $|C'\setminus \{w\}| = |C^*|$ and $C^*$ is a minimal vertex cover of $G$, it follows that $C' \setminus \{w\}$ is a minimal cover of $G$ since $G$ is well-covered.  Consequently,  $S' = V(G)\setminus (C'\setminus \{w\})$ is a maximal independent set of $G$.  On the other hand, since $y,x,v,z\notin S'$ and $N_G(y) = \{x,z\}$, it follows that $\{y\} \cup S'$ is an independent set of $G$, a contradiction. Thus, 
$(C',\delta')$ is minimal support,  as claimed.  

Together with Lemma \ref{cover},  this claim yields $\bht(I(G_\omega)) \geqslant |C'|=|C^*|+1 = \hgt(I(G_\omega)) + 1$. Thus, $I(G_\omega)$ is not unmixed,  a contradiction, and thus $m_1 \leqslant m$.  By the same way,  we get $m_1 \leqslant n$,  and $(1)$ follows.

$(2)$ From Part $(1)$ and our assumption, we have $m_k \leqslant \min\{m_1,\ldots,m_k,m,n\}$,  and hence
$$I(G_\omega) \colon y_k^{m_k}  = (y^pz^p, z^qu^q, u^rv^r, x^{m_k},  \ldots).$$

Since $N[y_k] \cap V(C) = \{x\}$ and $\deg_G(y) = \deg_G(v) = 2$,  we imply that the four monomials in the representation above are among minimal generators of $I(G_\omega) \colon y_k^{m_k}$ and the remaining minimal generators of $I(G_\omega) \colon y_k^{m_k}$ are not involving both $y$ and $v$.  Note that $y_kz, y_ku\notin E(G)$,  so $z^i,  u^i \notin I(G_\omega) \colon y_k^{m_k}$ for every $i$.  By Lemma \ref{CM-L01},  we obtain $q \leqslant \min \{p,r\}$,  as required.
\end{proof}

In the following lemmas, we use the setting as illustrated in Figure \ref{structureG}. Let $G_\omega$ be an unmixed edge-weighted graph where $G$ is in the class $\pc$ and let $C$ be a basic $5$-cycle of $G$.  Assume that $E(C) = \{xy,yz,zu,uv,vx\}$ and $\deg(x) > 2$.  Set
$$m = \omega(xy), p = \omega(yz), q=\omega(zu), r =\omega(uv), \text{ and } n = \omega(zx).$$

The aim of these lemmas is to show that each basic $5$-cycle of $G$ has a balanced vertex.

\begin{lem} \label{CM-L04} If $q < r$,  then $n\leqslant \min\{r,m\}$.
\end{lem}
\begin{proof}  Since $q \leqslant \min\{p,r\}$ by Lemma \ref{CM-L03},  we have
$$I(G_\omega)\colon u^q = (x^my^m, u^{r-q}v^r, v^nx^n, \ldots).$$

Since $q < r$,  by using the same argument as in the proof of Part $(2)$ of Lemma \ref{CM-L03} above,  we get $m\geqslant n$.  

We next prove that $n \leqslant r$.   Assume on the contrary that $n > r$.  Let $S$ be a maximal independent set of $G$ containing $x$ and $u$ and let $C^* = V(G)\setminus S$.  Then,  $C^*$ is a minimal vertex cover of $G$.  Let $\delta \colon C^* \to \Z_{>0}$ be a function such that $(C^*,\delta)$ is a minimal weighted vertex cover of $G_\omega$. Note that $x, u \notin C^*$ and $y, z, v \in C^*$.

Let $C' = C^* \cup\{u\}$ and $\delta' \colon C'\to \Z_{>0}$ given by 
$$\delta'(w) = \begin{cases} r & \text{ if } w = u,\\
n & \text{ if } w = v,\\
\delta(w) & \text{ otherwise}.
\end{cases}
$$

Then, $(C', \delta')$ is  a weighted vertex cover of  $G_\omega$. In fact, it suffices to check the condition of a weighted vertex cover of $G_\omega$ for the set of edges $\{zu, uv, xv\}$. 

If $e = zu$, then $z \in C'$ and $\delta'(z) = \delta (z) \leqslant \omega(zu)$ (the last inequality holds by look at the weighted vertex cover $(C^*, \delta)$).

If $e = uv$, then $u \in C'$ and $\delta'(u) = r = \omega(uv)$.

If $e = xv$, then $v \in C'$ and $\delta'(v) = n = \omega(xv)$. 

Next, we prove  $(C', \delta')$ is a minimal support weighted cover of $G_\omega$.  Assume on the contrary, there is a vertex $w \in C'$ such that  $(C'\setminus \{w\},\delta')$ is still a weighted vertex cover of $G_\omega$. Since $x \notin C'$ then $w$ could not be $y$ and $v$. We consider the following cases:

If $w = u$, then $\delta'(v) = n > r = \omega(uv)$, a contradiction. 

If $w = z$, then $\delta'(u) = r > q = \omega(zu)$, a contradiction. 

In other cases, since $|C'\setminus \{w\}| = |C^*|$,  $C^*$ is a minimal vertex cover of $G$, and $G$ is well-covered, then  $C' \setminus \{w\}$ is a minimal cover of $G$. Thus,  $S' = V(G)\setminus (C'\setminus \{w\})$ is a maximal independent set of $G$.  On the other hand, since $y,v,z,u\notin S'$ and $N_G(u) = \{v,z\}$, it follows that $\{u\} \cup S'$ is an independent set of $G$, a contradiction.

 Thus, $(C',\delta')$ is minimal support,  as claimed.

Since $|C'| = |C^*|+1$,  we have $\bht(I(G_\omega)) \geqslant |C^*|+1 = \hgt(I(G_\omega))+1$. This contradicts the fact that $I(G_\omega)$ is unmixed.  Therefore,  $n \leqslant r$,  as required.
\end{proof}

\begin{lem} \label{CM-L06} If $p=q<r$ and $n < m$,  then $p\leqslant m$.
\end{lem}
\begin{proof}  Assume on the contrary that $m < p$.  Let $S$ be a maximal independent set of $G$ containing $x$ and $z$ and let $C^* = V(G)\setminus S$ so that $C^*$ is a minimal vertex cover of $G$.  Let $\delta \colon C^* \to \Z_{>0}$ be a function such that $(C^*,\delta)$ is a minimal weighted cover of $G_\omega$.  

Let $C' = C^* \cup\{x\}$ and $\delta' \colon C'\to \Z_{>0}$ given by 
$$\delta'(w) = \begin{cases} p & \text{ if } w = y,\\
m & \text{ if } w = x,\\
\delta(w) & \text{ otherwise}.
\end{cases}
$$
By the same argumnet as the above Lemma, we can verify that $(C', \delta')$ is a minimal support weighted cover of $G_\omega$.  Since $|C'| = |C^*|+1$,  we have $\bht(I(G_\omega)) \geqslant |C^*|+1 = \hgt(I(G_\omega))+1$. This contradicts the fact that $I(G_\omega)$ is unmixed.  Therefore,  $p \leqslant m$,  as required.
\end{proof}

\begin{lem}\label{CM_L07} $C$ has a balanced vertex in the set $\{x,z,u\}$.
\end{lem}
\begin{proof} By Lemma \ref{CM-L03} we have $q \leqslant \min\{p,r\}$.  If $q < \min\{p,r\}$,  then $n=m$, $n \leqslant r$ and $m\leqslant p$ by Lemma \ref{CM-L04}.  Thus,  $x$ is a balanced vertex,  and thus  it remains to prove the lemma in the case $q =  \min\{p,r\}$.  We may assume that $p=q \leqslant r$.  We now consider two possible cases:

{\it Case $1$}: $p = q = r$.  By symmetry, we may assume that $m \leqslant n$.  We first claim that $\min\{m,n\} = m \leqslant p$.  Indeed,  assume on the contrary that $m > p$.  Let $S$ be a maximal independent set of $G$ containing $x$ and $u$ and let $C^* = V(G)\setminus S$.  Then,  $C^*$ is a minimal vertex cover of $G$.  Let $\delta \colon C^* \to \Z_{>0}$ be a function such that $(C^*,\delta)$ is a minimal vertex cover of $G_\omega$.

Let $C' = C^* \cup\{u\}$ and $\delta' \colon C'\to \Z_{>0}$ given by 
$$\delta'(w) = \begin{cases} m & \text{ if } w = y,\\
\min\{\delta(z),p\} & \text{ if } w = z,\\
p & \text{ if } w = u,\\
n & \text{ if } w = v,\\
\delta(w) & \text{ otherwise}.
\end{cases}
$$
By the same manner of the proof in the part $(1)$ of Lemma \ref{CM-L03},  $(C', \delta')$ is  a weighted vertex cover of  $G_\omega$.  It is straightforward to verify that it is a minimal support weighted cover of $G_\omega$.  Since $|C'| = |C^*|+1$,  we have $\bht(I(G_\omega)) \geqslant |C^*|+1 = \hgt(I(G_\omega))+1$. This contradicts the fact that $I(G_\omega)$ is unmixed.  Therefore,  $m \leqslant p$,  as claimed.

If $n\geqslant p$,  then $u$ is a balanced vertex on $C$. 
 
If $n < p$, we assume that $m < n$.  Let $S$ be a maximal independent set of $G$ containing $x$ and $u$ and let $C^* = V(G)\setminus S$.  Then,  $C^*$ is a minimal vertex cover of $G$.  Let $\delta \colon C^* \to \Z_{>0}$ be a function such that $(C^*,\delta)$ is a minimal weighted vertex cover of $G_\omega$.

Let $C' = C^* \cup\{x\}$ and $\delta' \colon C'\to \Z_{>0}$ given by 
$$\delta'(w) = \begin{cases} n & \text{ if } w = x,\\
p & \text{ if } w = v,\\
\delta(w) & \text{ otherwise}.
\end{cases}
$$
Again, $(C', \delta')$ is  a weighted vertex cover of $G_\omega$. Moreover, it is a minimal support weighted cover of $G_\omega$.  Since $|C'| = |C^*|+1$,  we have $\bht(I(G_\omega)) \geqslant |C^*|+1 = \hgt(I(G_\omega))+1$. This contradicts the fact that $I(G_\omega)$ is unmixed. Thus,  $m < n$ is impossible so that $m=n$. In this case, $x$ is a balanced vertex on $C$.

\medskip

{\it Case $2$}: $p = q < r$.  Then,  $n \leqslant \min\{r,m\}$ by Lemma \ref{CM-L04}.  In particular,  $n \leqslant m$.  If $n < m$,  then $p \leqslant m$ by Lemma \ref{CM-L06}, and then $z$ is a balanced vertex on $C$.  In the case $n=m$,  we have $z$ is  a balanced vertex on $C$ if $m \geqslant p$ or $x$ is a balanced vertex on $C$ if $m \leqslant p$.  The proof of the lemma is complete.
\end{proof}

\begin{lem}\label{CM_L08} Assume further that $\deg(z) > 2$. Then,  either $x$ or $z$ is a balanced vertex on $C$.
\end{lem}

\begin{proof} Since $\deg(x) > 2$ and $\deg(z) > 2$,  by Lemma \ref{CM_L07},  $C$ has a balanced vertex in the set $\{x,z,u, v\}$.  If $v$ is a balanced vertex, then
$$r = n \ \text { and } n \leqslant m \geqslant p \leqslant q \geqslant r.$$
On the other hand,  since $\deg(x) > 2$ and $\deg(z) > 2$,  by Lemma \ref{CM-L03} we obtain
$$p \geqslant q \leqslant r \text{ and }  m \geqslant n \leqslant r.$$

From those inequalities, we get $n=p=q=r\leqslant m$. Hence, $z$ is also a balanced vertex on $C$. In the same way, if $u$ is a balanced vertex, then $x$ is a balanced vertex on $C$ as well. Therefore, we conclude that either $x$ or $z$ is a balanced vertex on $C$.
\end{proof}

We are now in a position to prove the main result of the paper.

\begin{thm}\label{main-theorem} Let $G$ be a graph of girth at least $5$ and $\omega$ is a weight edge on $G$. Then, the following conditions are equivalent:
\begin{enumerate}
\item $G_\omega$ is Cohen-Macaulay.
\item $G$ is Cohen-Macaulay and $G_\omega$ is unmixed.
\item $G$ is in the class $\pc$ and the weight edge $\omega$ on $G$ satisfies:
\begin{enumerate}
\item The weight of any pendant edge in $G$ is greater than or equal to the weight of every edge adjacent to it.
\item Every basic $5$-cycle $C$ of $G$ has a balanced vertex adjacent to two vertices on $C$ of degree $2$.
\item If a vertex $x$ is on a basic $5$-cycle $C$ with $\deg_G(x)\geqslant 3$ and $N_C(x) = \{y,v\}$, then $\min\{\omega(xy),\omega(xv)\} \geqslant \max\{\omega(xw) \mid w \in N_G(x)\setminus\{y,v\}\}$.
\end{enumerate}
\end{enumerate}
\end{thm}

\begin{proof} $(1)\Longrightarrow (2)$ Since $G_\omega$ is Cohen-Macaulay, then $G_\omega$ is unmixed. On the other hand, $I(G)=\sqrt{I(G_\omega)}$,  by \cite[Theorem 2.6]{HTT}, $G$ is Cohen-Macaulay.

$(2)\Longrightarrow (3)$ Since $G$ is a Cohen-Macaulay graph of girth at least $5$, by Lemma \ref{HMT}, $G$ is in the class $\pc$. Now, we consider two following cases: First, in the case $G$ is just a $5$-cycle, we only need to prove the property $(b)$ and it follows immediately from \cite[Theorem 4.4]{PS}. Second, in the remain cases, i.e. $G$ is not a $5$-cycle, then the property $(a)$ equivalent to this statement: "For every pendant edge $xy$ of $G$ with $y$ is a leaf, then $\omega(xy) \geqslant \omega(xz)$ for any $xz \in E(G)$", and it follows from Lemma \ref{CM-L01}. In addition, the property $(b)$ follows from Lemma \ref{CM_L07}, and the property $(c)$ follows immediately from Lemma \ref{CM-L03}$(1)$.

$(3)\Longrightarrow (1)$ We prove by induction on the number of basic $5$-cycles of $G$. 

 If $G$ has no basic $5$-cycle, then its pendant edges form a perfect matching in $G$. In this case, combine with the condition $(a)$ and \cite[Lemma 5.3]{PS}, we get $G_\omega$ is Cohen-Macaulay.
 
 \medskip

 Assume that $G$ has  some basic $5$-cycles. If $G$ is just a $5$-cycle, by \cite[Theorem 4.4]{PS}, $G_{\omega}$ is Cohen-Macaulay as desired. If not then, assume $C_1,\ldots,C_r$ be the basic $5$-cycles of $G$ with $r\geqslant 1$ and $P$ be the set of pendant edges of $G$.  Assume that $E(C_1) = \{xy,yz,zu,uv,vx\}$ with 
$$\omega(xy) = m, \omega(yz) = p, \omega(zu) = q, \omega(uv) = r, \omega(vx) = n.$$

By our assumptions, $C_1$ has a balanced vertex such that two neighbors in $C_1$ are also of degree $2$. We may assume $x$ is such a vertex so that $m=n$ and $m\leqslant p\geqslant q \leqslant r \geqslant m$. Now we consider two possible cases:

\medskip

{\it Case $1$}: $\deg_G(x) = 2$. In this case, $N_G(x)=\{y,v\}$, and hence
$$I(G_\omega) \colon x^m = (y^m,  v^m, I((G_x)_\omega)) \ \text{ and } I(G_\omega)  +(x^m)=(x^m) + I((G\setminus x)_\omega).$$

Now, we will prove these ideals are Cohen-Macaulay. Observe that $G\setminus x$ is in the class $\pc$ with $r-1$ basic $5$-cycles $C_2,\ldots,C_r$ and pendant edges $P \cup \{zy,uv\}$ where $y$ and $v$ are leaves. We now verify the graph $(G\setminus x)_\omega$ satisfies the condition $(3)$. It suffices to prove the property $(a)$. Particularly, we only need to verify this property for the pendant edges $zy$ and $uv$. In particular, we will prove this property for the pendant edge $zy$, and similarly for the pendant edge $uv$. Let $zw \in E((G \setminus x)_\omega)$ for some $w \in V(G\setminus x) \setminus \{y\}$. If $w = u$, then by using condition of a balanced vertex $x$ in $C_1$, we have $\omega(zu) \leqslant \omega(zy)$. If $w \neq u$, then $w \notin C_1$. By applying Lemma \ref{CM-L03} on the basic $5$-cycle $C_1$, we get $\omega(zw) \leqslant \omega(zy)$. Thus, the property holds for the graph $(G\setminus x)_\omega$.  By the induction hypothesis,  $(G\setminus x)_\omega$ is Cohen-Macaulay,  so that $I(G_\omega) +(x^m)$ is Cohen-Macaulay.

In the same way,  we will prove that $I(G_\omega) \colon x^m$ is Cohen-Macaulay as follows. Since $C_1$ is a basic $5$-cycle then one of the vertices of $\{z,  u\}$ is a leaf in $G_x = G\setminus \{x,y,v\}$. Thus,  $G_x$ is in the class $\pc$ with $r-1$ basic $5$-cycles $C_2,\ldots,C_r$ and pendant edges $P \cup \{zu\}$. We now verify the graph $(G_ x)_\omega$ satisfies the condition $(3)$. It suffices to prove the property $(a)$. Particularly, it remains to verify this property for the pendant edge $zu$. If both vertices $z$ and $u$ are leaves, then nothing to do. Otherwise, assume $u$ is a leaf and $z$ is not. Let $zw$ be any edge in $E((G_x)_\omega)$, it follows that $w$ is not in the basic $5$-cycle $C_1$. Once again, applying Lemma \ref{CM-L03} on the basic $5$-cycle $C_1$, we get $\omega(zw) \leqslant \omega(zu)$. Thus, the property holds for the graph $(G_x)_\omega$, and hence $(G_x)_\omega$ is Cohen-Macaulay by the induction hypothesis.  Therefore, $I(G_\omega) \colon x^m$ is Cohen-Macaulay, too.

Since $\sqrt{I(G_\omega)  +(x^m)} = (x, I(G\setminus x))$ is Cohen-Macaulay,  it forces $G\setminus x$ is well-covered.  Since $x$ is not an isolated vertex,  it is a shedding vertex.   Moreover,

$\sqrt{I(G_\omega) \colon x^m} = (y,  v, I(G_x)) = I(G)\colon x$.  By Lemma \ref{dim}, we have $$\dim R/I(G_\omega) = \dim R/I(G_\omega)\colon x^m = \dim R/(I(G_\omega),x^m).$$ 
This implies that $I(G_\omega)$ is Cohen-Macaulay by Lemma \ref{CM-Q}.

\medskip

{\it Case $2$}: $\deg_G(x) > 2$.  Let $N(x) =\{y,v,y_1,\ldots,y_k\}$.  Since $m\geqslant m_i$ for all $i$ by Lemma \ref{CM-L03}, we obtain
$$I(G_\omega) \colon x^m = (y^m, v^m) +(y_1^{m_1}, \ldots, y_k^{m_k}, I(G\setminus \{x,y,v\})_\omega)$$
and
$$I(G_\omega)+( x^m) = (x^m, x^{m_1} y_1^{m_1}, \ldots, x^{m_k}y_k^{m_k}, I(G\setminus x)_\omega).$$

We now will prove these ideals are Cohen-Macaulay. Let $w$ be a new vertex and $H$ be a graph which is  obtained from $G$ by removing two edges $xy$ and $xv$ but adding a new edge $xw$. It means that $H$ is a graph with $V(H) = V(G) \cup \{w\}$ and $E(H) = (E(G) \cup \{xw\}) \setminus \{xy,xv\}$. Since $C_1$ is a basic $5$-cycle and $\deg_G(x) > 2$, then $\deg_G(y) = \deg_G(v) = 2$. Thus,  $w,y,v$ are leaves in $H$. Then $H$ is in the class $\pc$ with $r-1$ basic $5$-cycles and pendant edges $P \cup\{xw, uv\}$. Now we define the weight edge on $H$ by sending
$$e \mapsto \begin{cases} m & \text{ if } e = xw,\\
\omega(e) & \text{ otherwise},
\end{cases}
$$
which is still denoted by $\omega$.

 We now verify that $H_\omega$ satisfy the condition $(3)$. It suffices to prove the property $(a)$. In order to do this, it remains to verify this property for the pendant edges $xw$ and $uv$. It follows from Lemma \ref{CM-L03} (for the pendant edge $e = uv$) and the way we define the weight edge on $H$, $\omega(xw) = m$ (for the pendant edge $e = xw$). Thus, by the induction hypothesis,  $H_\omega$ is Cohen-Macaulay. Since $xw$ is an pendant edge of $H_\omega$, so that $x^mw^m \in I(H_\omega)$, by Lemma \ref{CM-Q} we have $I(H_\omega) \colon w^m$ is Cohen-Macaulay. Note that 
$$I(H_\omega)\colon w^m = (x^m, x^{m_1}y_1^{m_1}, \ldots, x^{m_k}y_k^{m_k}, I(G\setminus x)_{\omega})=I(G_\omega)+( x^m).$$
Hence,  $I(G_\omega)+( x^m)$ is Cohen-Macaulay.

In order to prove $I(G_\omega)\colon x^m$ is Cohen-Macaulay we use the same technique as above.  Let $H'$ be a graph with $V(H') = V(G\setminus \{y,v\})\cup\{w\}$ and $E(H') = E(G\setminus \{y,v\}) \cup \{xw\}$.  Next,  define the weight edge on $H'$, by sending
$$e \mapsto \begin{cases} m & \text{ if } e = xw,\\
\omega(e) & \text{ otherwise},
\end{cases}
$$
which is still denoted by $\omega$.

With this setting,  $H'_\omega$ is Cohen-Macaulay by the same argument as the previous case.  Thus,
$$I(H'_\omega) \colon x^m = (w^m, y_1^{m_1},\ldots, y_k^{m_k}, I((G\setminus\{x,y,v\})_\omega)$$
is Cohen-Macaulay by Lemma \ref{CM-Q}.  In particular,  $(y_1^{m_1},\ldots, y_k^{m_k}, I((G\setminus\{x,y,v\})_\omega)$ is Cohen-Macaulay,  and hence
$I(G_\omega) \colon x^m = (y^m,v^m) + (y_1^{m_1},\ldots, y_k^{m_k}, I((G\setminus\{x,y,v\})_\omega))$ is Cohen-Macaulay as well.

Finally, since $$\sqrt{I(G_\omega)\colon x^m} = (y,v, y_1,\ldots,y_k) + I(G\setminus\{x,y,v\}) = I(G)\colon x$$
and $$\sqrt{I(G_\omega)+(x^m)} = (I(G),x),$$
by the same argument as in Case $1$, we have $$\dim R/I(G_\omega) = \dim R/I(G_\omega)\colon x^m = \dim R/(I(G_\omega),x^m).$$ 
Therefore, $I(G_\omega)$ is Cohen-Macaulay by Lemma \ref{CM-Q},  and the proof is complete. 
\end{proof}

\begin{figure}[h]
\includegraphics[scale=0.45]{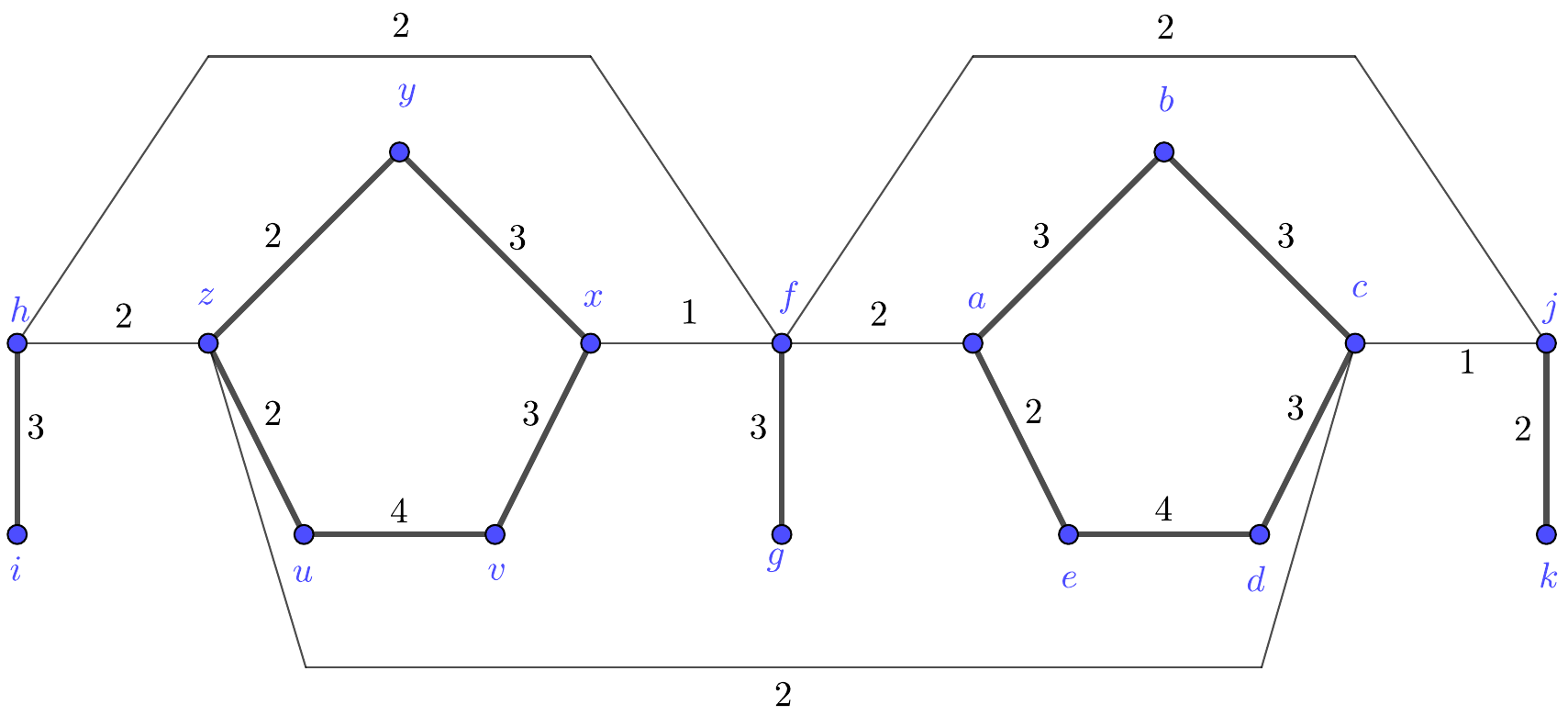}\\
\medskip
\caption{The Cohen-Macaulay edge-weighted graph.}
\label{CMWE}
\end{figure}

\begin{exm} The edge-weighted graph $G_\omega$ as depicted in Figure \ref{CMWE} is Cohen-Macaulay.
\end{exm}

 Indeed, we see from the figure that the underlying graph $G$ is in the class $\pc$ with three pendant edges $fg, hi, jk$; and two basic $5$ cycles $C_1:\ x\to y\to z\to u\to v\to x$ and $C_2: \ a\to b\to c\to d\to e\to a$. Note that $z$ is a balanced vertex on $C_1$ and $c$ is the one on $C_2$ they satisfy the condition $(b)$ in Theorem \ref{main-theorem}.

We can easily verify that the conditions $(a)-(c)$ in Theorem \ref{main-theorem} holds for $G_\omega$, and thus $G_\omega$ is Cohen-Macaulay.

\medskip

% -----------------------------------------------------------
\subsection*{Acknowledgment}  This work is  partially supported by NAFOSTED (Vietnam) under the grant number 101.04-2023.36.

% -----------------------------------------------------------


\begin{thebibliography}{FHN1}
\bibitem {BC} T. Biyik\u{o}glu and Y. Civan, {\it Vertex-decomposable graphs, codismantlability, Cohen-Macaulayness, and Castelnuovo-Mumford regularity},  Electron. J. Combin. {\bf 21} (1) (2014) Paper 1.1, 17 pp.

% \bibitem {CHHKTT} G.  Caviglia,  H. T.  Ha,  J.  Herzog,  M.  Kummini,  N.  Terai and N.V.  Trung,  {\it Depth and regularity modulo a principal ideal},  J. Algebraic Combin. {\bf 49} (2019), 1-20.

\bibitem {MDV} L.T.K. Diem, N.C. Minh and T. Vu, {\it The sequentially Cohen-Macaulay property of edge ideals of edge-weighted graphs}, arXiv:2308.05020.

\bibitem {DT} L.X. Dung and T.N.  Trung,  {\it Cohen-Macaulay oriented graphs with large girth}, arXiv:2308.11907.

%\bibitem {HHZ} J. Herzogs, T. Hibi and X. Zheng, \textit{Cohen-Macaulay chordal graphs}, J. Combin. Theory Ser. A, 113 (2006), no. 5, 911-916.

\bibitem {HLMRV} H.T. Ha, K. Lin, S. Morey, E. Reyes and R. H. Villarreal, {\it Edge ideals of oriented graphs}, Int. J. Algebra Comput. {\bf 29} (2019), 535–559.

\bibitem {HTT} J. Herzog, Y. Takayama and N. Terai, {\it On the radical of a monomial ideal}, Arch. Math. {\bf 85} (2005), 397-408.

\bibitem {HMT} D.T.  Hoang,  N.C.  Minh and T.N. Trung,  {\it Cohen-Macaulay graphs with large girth},  J. Algebra Appl. {\bf 14} (2015),  no. 7,  1550112,  16 pp.

\bibitem {PS} C. Paulsen and S. Sather-Wagstaff,  {\it Edge ideals of weighted graphs},  J. Algebra Appl. {\bf 12} (2013), no 5, 1250223, 24pp.

\bibitem {PRT} Y. Pitones, E. Reyes, and J. Toledo, {\it Monomial ideals of weighted oriented graphs}, Electron. J. Combin., 26 (2019), no. 3, Research Paper P3.44.

\bibitem {PRV} Y. Pitones, E. Reyes and R. H. Villarreal, {\it Unmixed and Cohen-Macaulay weighted oriented K\"{o}nig graphs}, Studia Sci. Math. Hungar. {\bf 58} (2021), no. 3, 276-292.

\bibitem {P} M. D. Plummer, {\it Some covering concepts in graphs}, Journal of Combinatorial Theory, {\bf 8} (1970), 91-98.

\bibitem {SSTY} A.A.  Seyed Fakhari,  K.  Shibata, K.,  N.Terai,  and S.  Yassemi,  {\it Cohen–Macaulay edge-weighted edge ideals of very well-covered graphs},   Commun. Algebra {\bf 49}(10), 4249-4257 (2021).

\bibitem {Vi} R. Villarreal, {\it Monomial Algebras}, Monographs and Textbooks in Pure and Applied Mathematics Vol. 238, Marcel Dekker, New York, 2001.

\bibitem {W} R. Woodroofe, {\it Vertex decomposable graphs and obstructions to shellability}, Proc. Amer. Math. Soc. {\bf 137} (2009), no. 10, 3235-3246.
\end{thebibliography}
\end{document}